\documentclass[10pt]{article}
\usepackage{lmodern}
\usepackage[T1]{fontenc}
\usepackage{amsmath}
\usepackage{amsthm}
\usepackage{amssymb}
\usepackage{titlefoot}
\usepackage[small]{titlesec}
\usepackage{listings}
\usepackage{authblk}

\DeclareFixedFont{\ttb}{T1}{txtt}{bx}{n}{9.5} % for bold
\DeclareFixedFont{\ttm}{T1}{txtt}{m}{n}{10}  % for normal

\usepackage{color}
\definecolor{deepblue}{rgb}{0,0,0.5}
\definecolor{deepred}{rgb}{0.6,0,0}
\definecolor{deepgreen}{rgb}{0,0.5,0}

\newcounter{todocounter}

\newcommand\pythonstyle{\lstset{
	tabsize=1,
	language=Python,
	basicstyle=\ttm,
	otherkeywords={sage},             % Add keywords here
	keywordstyle=\ttb\color{deepblue},
	emph={True},          % Custom highlighting
	emphstyle=\ttb\color{deepred},    % Custom highlighting style
	stringstyle=\color{deepgreen},
	frame=tb,                         % Any extra options here
	showstringspaces=false            %
}}

\lstnewenvironment{python}[1][]{
	\pythonstyle
	\lstset{#1}
}{}

\newcommand\pythoninline[1]{{\pythonstyle\lstinline!#1!}}

\theoremstyle{plain}
\newtheorem{theorem}{Theorem}[section]

\newtheorem{conjecture}[theorem]{Conjecture}

\makeatletter
\newfont{\footsc}{cmcsc10 at 8truept}
\newfont{\footbf}{cmbx10 at 8truept}
\newfont{\footrm}{cmr10 at 10truept}
\renewenvironment{abstract}%
                {
                  \begin{list}{}%
                     {\setlength{\rightmargin}{1in}%
                      \setlength{\leftmargin}{1in}}%
                   \item[]\ignorespaces\begin{small}}%
                 {\end{small}\unskip\end{list}}

\newcommand{\eval}[2][\right]{\relax\ifx#1\right\relax \left.\fi#2#1\rvert}

\datefoot{\today}
\amssubj{05C60, 68R15}

\newpagestyle{main}[\small]{
        \headrule
        \sethead[\usepage][][]
        {\sc On Isomorphism Classes of Generalized Fibonacci Cubes}{}{\usepage}}

\title{\sc On Isomorphism Classes of Generalized Fibonacci Cubes}
\author[1]{Jernej Azarija}
\author[1,2,3]{Sandi Klav\v{z}ar}
\author[4]{Jaehun Lee}
\author[5]{Jay Pantone}
\author[4]{Yoomi Rho}
\affil[1]{Institute of Mathematics, Physics, and Mechanics, Ljubljana, Slovenia\hspace*{-1cm}\newline\emph{jernej.azarija@gmail.com}\vspace*{.3cm}}
\affil[2]{Faculty of Mathematics and Physics, University of Ljubljana, Slovenia\hspace*{-1cm}}
\affil[3]{Faculty of Natural Sciences and Mathematics, University of Maribor, Slovenia\hspace*{-.5cm}\newline\emph{sandi.klavzar@fmf.uni-lj.si}\vspace*{.3cm}}
\affil[4]{Department of Mathematics, Incheon National University, Korea\hspace*{-1.5cm}\newline\emph{dlwogen257@naver.com}, \emph{rho@incheon.ac.kr}\vspace*{.3cm}}
\affil[5]{Department of Mathematics, University of Florida, USA\hspace*{-2cm}\newline\emph{jaypantone@ufl.edu}}

\titleformat{\section}
        {\large\sc}
        {\thesection.}{1em}{}

\date{}

\begin{document}
\maketitle

\pagestyle{main}

\begin{abstract}
The generalized Fibonacci cube $Q_d(f)$ is the subgraph of the $d$-cube $Q_d$ induced on the set of all strings of length $d$ that do not contain $f$ as a substring.  It is proved that if $Q_d(f) \cong Q_d(f')$ then $|f|=|f'|$. The key tool to prove this result is a result of Guibas and Odlyzko about the autocorrelation polynomial associated to a binary string. It is also proved that there exist pairs of strings $f, f'$ such that $Q_d(f) \cong Q_d(f')$, where $|f| \ge \frac{2}{3}(d+1)$ and $f'$ cannot be obtained from $f$ by its reversal or binary complementation. Strings $f$ and $f'$ with $|f|=|f'|=d-1$ for which $Q_d(f) \cong Q_d(f')$ are characterized. 
\end{abstract}

%%%%%%%%%%%%%%%%%%%%%%%%%%%%%%%%%%%%%%%%%%%%%%%%%%%%%%%
\section{Introduction}
\label{sec:intro}
%%%%%%%%%%%%%%%%%%%%%%%%%%%%%%%%%%%%%%%%%%%%%%%%%%%%%%%

An element of $\{0,1\}^d$ is called a {\em binary string} (henceforth just called a \emph{string}) of length $d$, with the usual concatenation notation. For example, $0^{d-1}1$ is the string of length $d$ consisting of $d-1$ 0 bits followed by a single 1 bit. We will denote by $e_i = 0^{i-1}10^{d-i}$ the $i^{\rm th}$ unit string in $\{0,1\}^d$. 

Let $d\ge 1$ be a fixed integer. The \emph{$d$-cube} $Q_d$ is the graph whose vertices are the binary strings of length $d$, with an edge connecting vertices $v_1$ and $v_2$ if the underlying strings differ in exactly one position. Given a graph $G$, the set of vertices of $G$ is denoted by $V(G)$. We use $d_G(u,v)$ to denote the length of the shortest path connecting $u$ and $v$ in $G$. Lastly, we will write $G\cong H$ to signify that the graphs $G$ and $H$ are isomorphic. 

For a given string $f$ and integer $d$, the \emph{generalized Fibonacci cube} $Q_d(f)$ is the subgraph of $Q_d$ induced by the set of all strings of length $d$ that do not contain $f$ as a consecutive substring. Indeed, this generalizes the notion of the $d$-dimensional \emph{Fibonacci cube} $\Gamma_d = Q_d(11)$, which is the graph obtained from the $d$-cube $Q_d$ by removing all vertices that contain the substring $11$.

Fibonacci cubes were introduced by Hsu~\cite{hsu-93} as a model for interconnection networks. Like the hypercube graphs, Fibonacci cubes have several properties which make them ideal as a network topology, yet their size grows significantly slower than that of the hypercubes. Fibonacci cubes have been extensively investigated; see, for example, the recent survey by Klav\v{z}ar~\cite{kl-2013} and even more recent papers of Klav\v{z}ar and Mollard~\cite{klmo-2014} and Vesel~\cite{ve-2014}. In the first of these papers, different asymptotic properties of Fibonacci cubes are established, while in the latter a linear recognition algorithm is designed for recognizing Fibonacci cubes, improving the previous best recognition algorithm of Taranenko and Vesel~\cite{tave-2007}.

Later, Ili\'c, Klav\v{z}ar, and Rho~\cite{ilkl-2012a} introduced the idea of generalized Fibonacci cubes (as defined above). Under the same name, the graphs $Q_d(1^s)$ were studied by Liu, Hsu, and Chung~\cite{lihs-1994} and Zagaglia Salvi~\cite{za-1996}. The analysis of the properties of generalized Fibonacci cubes led to the study of several problems related to the combinatorics of words. To study their isometric embeddability into hypercubes, good and bad words were introduced by Klav\v{z}ar and Shpectorov~\cite{klsh-2012}, where it was proved that about eight percent of all words are good. Isometric embeddability and hamiltonicity of generalized Fibonacci cubes motivated the ideas of the index and parity of a binary word, as defined by Ili\'c, Klav\v{z}ar, and Rho~\cite{ilkl-2012b, ilkl-2012c}.

In this paper we consider the following fundamental question about the generalized Fibonacci cubes: for which binary strings $f$ and $f'$ and positive integers $d$ are the generalized Fibonacci cubes $Q_d(f)$ and $Q_d(f')$ isomorphic? It is easy to see that if $f'$ is the binary complement of $f$, or if $f'$ is the reverse of $f$ (the reverse of $f=f_1\ldots f_d$ is $f_d\ldots f_1$), then $Q_d(f) \cong Q_d(f')$ for any dimension $d$. Hence we say that a  pair of binary strings $f, f'$ is {\em trivial} if $f'$ can be obtained from $f$ by binary complementation, reversal, or composition of these mappings. We are therefore only interested in the behavior of the other pairs, which we call the {\em non-trivial pairs}.

We proceed as follows. In the next section, we prove that if $Q_d(f) \cong Q_d(f')$, then $|f|=|f'|$. We also prove that there exist non-trivial pairs of strings $f, f'$ such that $Q_d(f) \cong Q_d(f')$, where $|f| \ge \frac{2}{3}(d+1)$. In the last section, we prove that if $|f| = d-1$, then $Q_d(f) \cong Q_d(f')$ if and only if $f$ and $f'$ have the same block structure. Several conjectures are posed along the way.

%%%%%%%%%%%%%%%%%%%%%%%%%%%%%%%%%%%%%%%%%%%%%%%%%%%%%%%
\section{The Length of Forbidden Words}
\label{sec:length}
%%%%%%%%%%%%%%%%%%%%%%%%%%%%%%%%%%%%%%%%%%%%%%%%%%%%%%%

In this section we first prove that if $Q_d(f) \cong Q_d(f')$, then $|f|=|f'|$. Then we pose the question whether there is some relation between $|f|$ ($= |f'|$) and $d$ provided that $Q_d(f) \cong Q_d(f')$. To this end we prove that there exist non-trivial pairs $f, f'$ such that $Q_d(f) \cong Q_d(f')$ and $|f| \ge \frac{2}{3}(d+1)$. We also conjecture that for any non-trivial pair $f,f '$ such that $Q_d(f) \cong Q_d(f')$, we must have $|f| \ge \frac{2}{3}(d+1)$.

The {\em autocorrelation polynomial} $p_f(z)$ associated to a binary string $f = f_1\ldots f_k\in \{0,1\}^k$ is defined as 
$$p_f(z)= \sum_{i=0}^{k-1} c_i z^i\,,$$
where $c_i = 1$ if the length $k-i$ suffix of $f$ is equal to the length $k-i$ prefix of $f$, i.e., if $f_{i+1}\ldots f_{k} = f_1\ldots f_{k-i}$, and $c_i=0$ otherwise; see Flajolet and Sedgewick~\cite[p.~60]{Sedgewick}. Note that the autocorrelation polynomial of $f \in \{0,1\}^k$ is of degree at most $k-1$ and has degree $k-1$ if and only if the last bit of $f$ is equal to the first bit of $f$.  Observe also that
\begin{equation}
\label{eq:extremes}
p_{0^k}(z) = \sum_{i=0}^{k-1} z^i \quad \mbox{and} \quad p_{0^{k-1}1}(z) = 1\,.
\end{equation}
We note in passing that more generally, $p_{f}(z) = 1 + z + \cdots + z^{k-1}$ if and only if $f=0^k$ or $f=1^k$ and that $p_{f}(z) = 1$ if and only if every non-trivial suffix of $f$ is different from the prefix of $f$ of the same length. Such words were named {\em prime} in~\cite{ilkl-2012c}. 

The following theorem establishes that only strings which have the same length can generate isomorphic generalized Fibonacci cubes.
\begin{theorem}
\label{thm:equal-length}
If $|f|, |f'| \le d$ and $Q_d(f) \cong Q_d(f')$, then $|f|=|f'|$.
\end{theorem}

\proof
Set $n_d(f) = |V(Q_d(f))|$ and assume that $|f| < |f'|$. We will show that $n_d(f) < n_d(f')$, from which the theorem follows.

The key idea of the proof is to apply a result of Guibas and Odlyzko~\cite[p.~204]{Guibas} stating that if $g$ and $g'$ are binary strings such that $p_{g}(2) > p_{g'}(2)$, then $n_d(g) \geq n_d(g')$.  Moreover, since the values $n_d(g)$ only depends on $p_g(z)$, cf.~\cite[Proposition~1.4]{Sedgewick}, we also infer that if $|g| = |g'|$ and $p_{g}(2) = p_{g'}(2)$, then $n_d(g) \geq n_d(g')$. It follows that if $f$ is a binary string of length $k \leq d$, then
\begin{equation}
\label{eq:guibas}
n_d(0^{k-1}1) \leq n_d(f) \leq n_d(0^k)\,.
\end{equation}
In addition, since $0^{k}$ is a strict substring of $0^k1$ we infer that
\begin{equation}
\label{eq:just-some}
n_d(0^k) < n_d(0^k 1)\,.
\end{equation}
Combining~\eqref{eq:guibas} and~\eqref{eq:just-some} we conclude that
$$n_d(f) \leq n_d(0^k) <  n_d(0^k1) \leq n_d(f')\,,$$
where the last inequality is due to the assumption that $|f'|>k$.
\qed

% Let us remark that the validity of a more general equality over alphabets of arbitrary length is still an open problem.

Hence, non-trivial pairs $f, f'$ such that $Q_d(f) \cong Q_d(f')$ are of the same length. We next ask what is the relation of this length with the dimension of the corresponding generalized Fibonacci cubes. The following result could be the extremal case.

\begin{theorem}
\label{thm:3k-1}
If $k\ge 2$, then $Q_{d}(0^{k}1^{k}) \cong Q_{d}(0^{k+1}1^{k-1})$ for any $d\le 3k-1$.
\end{theorem}

\proof
Let $k\ge 2$ be a fixed integer and set $f = 0^{k}1^{k}$, $f' = 0^{k+1}1^{k-1}$, $G = Q_{3k-1}(f)$, and $G' =  Q_{3k-1}(f')$. Let $X = V(Q_{3k-1})\setminus V(G)$ and $X' = V(Q_{3k-1})\setminus V(G')$. For any $0 \leq i \leq k$ let
$$X_i = \{ufv:\ |u|=i, |v| = k-1-i \}\,.$$
Then, by definition,
	\[X = \bigcup_{i=0}^{k-1} X_i.\]
Let $w=w_1\ldots w_{3k-1}$ be an arbitrary vertex whose underlying string is in $X_i$. It then follows that $w_{k+1}\ldots w_{2k-1} = 0^{i}1^{k-i-1}$, and so $X_j\cap X_{k} = \emptyset$ holds for all $j\ne k$. Since $|X_i| = 2^{k-1}$, we have $|X| = k2^{k-1}$. With a parallel argument we infer that also $|X'| = k2^{k-1}$. This implies that $|V(G)| = |V(G')|$.

Consider now the mapping $\alpha: V(Q_{3k-1}) \rightarrow V(Q_{3k-1})$ defined by
$$\alpha(u_1\ldots u_{3k-1}) = u_1\ldots u_k \overline{u_{2k}} u_{k+1}\ldots u_{2k-1}u_{2k+1}\ldots u_{3k-1}\,.$$
In particular, $\alpha$ fixes the first $k$ and the last $k-1$ coordinates. Since transposition of coordinates, complementation of a coordinate, and any composition of such mappings are all automorphisms of a hypercube, $\alpha$ is an automorphism of $Q_{3k-1}$. Consider now $G$ and $G'$ as subgraphs of $Q_{3k-1}$ and the restriction $\alpha|_G$ of $\alpha$ to $G$. Since $|V(G)| = |V(G')|$, it remains to prove that $\alpha|_G: V(G) \rightarrow V(G')$.

Suppose on the contrary that for some $u\in V(G)$ we have $\alpha|_G(u) = w\notin V(G')$. Then $w=x0^{k+1}1^{k-1}y$, where $|x| = i$ for some $0 \leq i \leq k$, and  $|y| = k-i-1$.
% Note that $\alpha|_G^{-1}$ maps $w_{2k}$ to $\overline{w_{k+1}}$.
It is straightforward to see that $\alpha|_G^{-1}(w) = x0^{k}1^{k}y$. Since this is not a  vertex of $V(G)$ we have a contradiction.

We have thus proved the result for $d=3k-1$. Note that all of the above arguments also work for $2k+1\le d\le 3k-2$ and when $d \leq 2k$, the assertion is trivial.
\qed

Motivated by the last theorem we pose:

\begin{conjecture}
\label{conj:dim-minus-1}
If $f$ and $f'$ are binary strings such that $Q_d({f})\cong Q_d({f'})$, then $Q_{d-1}(f)\cong Q_{d-1}(f')$.
\end{conjecture}

\begin{conjecture}
\label{conj:2/3}
Let $f, f'$ be a non-trivial pair such that $Q_d(f) \cong Q_d(f')$. Then $|f| \ge \frac{2}{3}(d+1)$.
\end{conjecture}

We have verified these conjectures for small strings using the Sage package~\cite{Sage}. More precisely, we tested both conjectures for all $d\le 12$ and all non-trivial pairs $f, f'$. See Appendix~\ref{appendix} for the corresponding procedures.

%%%%%%%%%%%%%%%%%%%%%%%%%%%%%%%%%%%%%%%%%%%%%%%%%%%%%%%
\section{The Number of Blocks in Forbidden Words}
\label{sec:blocks}
%%%%%%%%%%%%%%%%%%%%%%%%%%%%%%%%%%%%%%%%%%%%%%%%%%%%%%%

In this section we characterize the binary strings $f,f'$ of length $d-1$ for which $Q_{d}(f)\cong Q_{d}(f')$. It turns out that they are precisely the strings with the same block structure. 

Let $\nu(f)$ denote one less than the number of blocks of $f=f_1f_2\ldots f_{|f|}$. For example $\nu(0110)=2$. When a bit is different from the previous bit we call its index an {\it index of bit change} and denote it by $i_j$. Therefore $f_{i_{j}-1}\neq f_{i_j}$ for $i_1<i_2<\cdots <i_{\nu(f)}$.

\begin{theorem}
\label{thm:d-1}
Let $d\ge 2$ and let $f$, $f'$ be binary strings of length $d-1$. Then $\nu({f})=\nu({f'})$ if and only if $Q_{d}(f)\cong Q_{d}(f')$.
\end{theorem}

\proof
Set $A=V(Q_d)\setminus V(Q_d(f))$.
Then $A=\{\overline{f_1}f, f_1f, ff_{|f|}, f\overline{f_{|f|}}\}$.
Similarly set $A'=V(Q_d)\setminus V(Q_d(f'))=\{\overline{f'_1}f', f'_1f', f'f'_{|f|}, f'\overline{f'_{|f'|}}\}$.
Let $2\le i_1<\cdots <i_{\nu(f)}\le d-1$ be the indices of bit change of $f$ and let $2\le i'_1<\cdots <i'_{\nu(f')}\le d-1$ be the indices of bit change of $f'$.
Since $(\overline{f_1}f)_{\tau}=(f_1f)_{\tau}=f_{\tau-1}$ for $2\le \tau\le d-1$ and $(ff_{|f|})_{\tau}=(f\overline{f_{|f|}})_{\tau}=f_{\tau}$ for $2\le \tau\le d-1$, the strings $f_1f$ and $ff_{|f|}$ are different precisely for $\tau=i_j$, $1\le j\le \nu(f)$. % and the same for $\theta\neq 1,d,i_j$ for $1\le j\le \nu(f)$.
It follows that $d_{Q_d}(f_1f, ff_{|f|})=\nu(f)$ and by a parallel argument $d_{Q_d}(f'_1f', f'f'_{|f'|})=\nu(f')$.

Assume first that $\nu({f})=\nu({f'})$. We may without loss of generality assume that $f_1=f'_1=0$. Then $f_{|f|}=f'_{|f|}$. Let $\phi$ be a permutation of $\{1,\ldots,d\}$ such that $\phi(i_j)=i'_j$, $\phi(1)=1$, $\phi(d)=d$. Set $\psi(x_1\ldots x_{d})=y_{1}\ldots y_{d}$, where
$$y_{\tau}=\left \{\begin{array}{ll}
x_{\phi^{-1}(\tau)}; & \text{if \ $f_{\phi^{-1}(\tau)}=f'_{\tau}$,} \\[4pt]
 {\overline{x_{\phi^{-1}(\tau)}}}; & \text{otherwise.}\end{array}\right.$$
As transposition of coordinates, complementation of a coordinate, and compositions of such mappings are automorphisms of a hypercube, $\psi$ is an automorphism of $Q_d$. Also $\psi$ sends the vertices of $A$ to $A'$. Thus  $\psi$ is an isomorphism from $Q_d(f)$ to $Q_d(f')$.

To prove the converse assume that $Q_{d}(f)\cong Q_{d}(f')$. We may without loss of generality assume that $\nu(f)\le \nu(f')$.

Assume $\nu(f)=0$. Then $f_1f = ff_{|f|}$ and hence $|A| = 3$, while $|A'|=4$ if $\nu(f')\neq 0$. Therefore $\nu(f')=0$.

Assume $\nu(f)=1$. Then note that the vertices from $A$ induce a path on four vertices  and hence $|E(Q_d(f))| = d2^{d-1}-(4d-3)$. If $\nu(f')> 1$, then the vertices from $A'$ induce two disjoint copies of $K_2$ and hence  $|E(Q_d(f'))| = d2^{d-1}-(4d-2) \ne |E(Q_d(f))|$. We conclude that $\nu(f')=1$.

For the rest of the proof we can thus assume that $\nu(f)\ge 2$. The subgraph of $Q_d$ induced on $A$ consists of two edges $\{\overline{f_1}f,  f_1f\}$ and $\{ff_{|f|}, f\overline{f_{|f|}}\}$ where $d( f_1f, ff_{|f|})= \nu(f)$ and $d(\overline{f_1}f, f\overline{f_{|f|}})= \nu(f)+2$. 
   Denote $\overline{f_1}f$, $f_1f$, $ff_{|f|}$, $f\overline{f_{|f|}}$ by $a$, $b$, $c$, $d$, respectively.
   Consider the shortest $b,c$-path constructed by changing from left to right the bits in which $b$ and $c$ differ: 
$$b=f_1f \rightarrow f_1f+e_{i_1} \rightarrow f_1f+e_{i_1}+e_{i_2} \rightarrow \cdots \rightarrow f_1f+e_{i_1}+\cdots+e_{i_{\nu(f)}}=ff_{|f|}=c,$$
where addition is taken modulo 2. 

   Denote the $j$-th internal vertex $f_1f+e_{i_1}+\cdots+e_{i_j}$ of this path by $x_j$ for $1\le j\le \nu(f)-1$.
   Similarly, denote $\overline{f'_1}f'$, $f'_1f'$, $f'f'_{|f|}$, $f'\overline{f'_{|f|}}$ by $a'$, $b'$, $c'$, $d'$, respectively.
   Set $k = \nu(f)-1$, $\ell = \nu(f')-1$, and recall that $k\ge 1$.
Let $\psi:Q_{d}(f)\rightarrow Q_{d}(f')$ be an isomorphism and let $x'_j = \psi(x_j)$ for $1\le j\le k$.

 %and set $x'_1, x'_2,\ldots, x'_{\ell}$ to be the internal vertices on the shortest path from $b'$ to $c'$, $(\overline{f'_1}f')(f'_1f')(f'_1f'+e_{i'_1})(f'_1f'+e_{i'_1}+e_{i'_2})\cdots ((f'_1f'+e_{i'_1}+\cdots+e_{i'_\nu(f')})=(f'f'_{|f|}))(f'\overline{f'_{|f|}})$ where $x'_i$ is closer to $b'$ than $x'_{i+1}$ is.
%Let $\psi(y_i)=y'_i$ for all $1\le i\le \ell$.

%From now on suppose $\nu(f)<\nu(f')$. We will prove the theorem by showing that we get a contradiction in each case.

Assume $k= 1$. Then $x_1$ is of degree $d-2$ and hence ${\rm deg}_{Q_d(f')}(x_1')=d-2$.
This means that $x'_1$ is adjacent to two vertices among $a',b',c',d'$.
As $Q_d$ is bipartite, $x'_1$ is adjacent to one of $a',b'$ and one of $c',d'$.
If $x'_1$ is adjacent to $b'$ and $c'$, then $\ell+1=d_{Q_d}(b',c')\le 2$ and hence $\ell=1$.
If $x'_1$ is adjacent to $a'$ and $c'$, then $\ell+2=d_{Q_d}(a',c')\le 2$, a contradiction.
Similarly we get contradictions if $x'_1$ is adjacent to $d'$.

Assume $k\ge 2$. Then the vertices $x_1$ and $x_k$ of $Q_d(f)$ are of degree $d-1$. Considering that $x_1 \rightarrow \cdots \rightarrow x_k$ is a path in $Q_d(f)$,  we see that $d_{Q_d(f)}(x_1,x_k)=k-1$. Therefore the vertices $x'_1$ and $x'_k$ of $Q_d(f')$ are of degree $d-1$ and $d_{Q_d(f')}(x'_1,x'_k)=k-1$. We distinguish three cases.

\medskip\noindent
{\bf Case 1}: ($x'_1$ and $x'_k$ are adjacent to a common vertex among $a',b',c',d'$)\\
Now, $k-1=d_{Q_d}(x'_1, x'_k)\le 2$ and hence $k\le 3$. Also, $k$ is odd as $Q_d$ is bipartite, and thus $k=3$.
Considering that $x_1x_2 x_3$ is a shortest path in $Q_d$, it follows that $x_1$ and $x_3$ have distance two and hence have a common neighbor which is different from $x_2$ in $Q_d$. Call it $u$. Then, $u$ is not $a$, $b$, $c$, or $d$ because of its distances from $x_1$ and $x_3$. Therefore $u\in Q_d(f)$. Set $u'=\psi(u)$. Then $u'\in Q_d(f')$ and hence $d_{Q_d}(x'_1, u')\le d_{Q_d(f')}(x'_1, u')=d_{Q_d(f)}(x_1, u)=1$, which means that $d_{Q_d}(x'_1, u') = 1$. Similarly, $d_{Q_d}(x'_3, u')=1$. Hence in $Q_d$, $x'_1$ and $x'_3$ have three common neighbors: $u'$, $x_2'$, and one of $a',b',c',d'$. This is a contradiction, because hypercubes are $K_{2,3}$-free.
\ \\\ \\
From now on we regard that $x'_1$ and $x'_k$ are not adjacent to a common vertex among $a',b',c',d'$.

\medskip\noindent
{\bf Case 2}: ($x'_1$ and $x'_k$ are either adjacent to $a'$ and $b'$ or adjacent to $c'$ and $d'$)\\
We may without loss of generality assume the first. Then $k-1=d_{Q_d}(x'_1, x'_k)\le 3$ and hence $k\le 4$. Also, $k$ is even as $Q_d$ is bipartite. Therefore $k=2$ or $k=4$. We distinguish two subcases.

\medskip\noindent
{\bf Case 2a}: ($k=2$)\\
It is well-known that in $Q_d$ a given edge lies in $d-1$ cycles of length 4.
Among the $4$-cycles containing the edge $x_1x_2$,
one contains $b$ and another contains $c$. This means that there are $d-3$ cycles of length 4 containing the edge $x_1x_2$ in $Q_d(f)$. Among the $4$-cycles containing the edge $x'_1x'_2$, only one contains $a'$ and $b'$ together, and no other contains $a',b',c'$, or $d'$. This means that there are $d-2$ cycles of length 4 containing the edge $x'_1x'_2$ in $Q_d(f')$, a contradiction.

\medskip\noindent
{\bf Case 2b}: ($k=4$)\\
It is known that for two given vertices at distance three, there are exactly three internally vertex-disjoint shortest paths in $Q_d$, and therefore there are such paths between $x_1$ and $x_4$. Let $R=x_1uvx_4$ be any one of them which is different from $x_1x_2x_3x_4$. Considering the distances of $u,v$ from $x_1, x_k$, we obtain that $u,v\in Q_d(f)$ and hence $R$ is a path in $Q_d(f)$. Therefore $\psi(R)$ is a path in $Q_d(f')$. By the assumption that $k=4$, there is also an $x'_1,x'_4$-path through $a'$ and $b'$, implying that there are (at least) four internally disjoint shortest paths between $x'_1$ and $x'_4$ in $Q_d$, which is a contradiction.

\medskip\noindent
{\bf Case 3}: ($x'_1$ is adjacent to one of $a'$ and $b'$ while $x'_k$ is adjacent to one of $c'$ and $d'$, or vice versa)

Firstly, assume that $x'_1$ is adjacent to $b'$ while $x'_k$ is adjacent to $c'$. Then 
\begin{eqnarray*}
\ell+1 & = & d_{Q_d}(b',c')\\
& \le & d_{Q_d}(b',x'_1)+d_{Q_d}(x'_1, x'_k)+ d_{Q_d}(x'_k,c') \\ 
& = & 2+d_{Q_d}(x'_1, x'_k)\\
& \le & 2+d_{Q_d(f')}(x'_1, x'_k) \\
& =& 2+d_{Q_d(f)}(x_1, x_k)\\
& = &k+1\,.
\end{eqnarray*}
Hence, under this assumption $\ell=k$. 

Alternatively, assume that $x'_1$ is adjacent to $a'$ while $x'_k$ is adjacent to $c'$. Then,
\begin{eqnarray*}
\ell+2 & = & d_{Q_d}(a',c')\\
& \le & d_{Q_d}(a',x'_1)+d_{Q_d}(x'_1, x'_k)+ d_{Q_d}(x'_k,c') \\ 
& \le & 2+d_{Q_d(f')}(x'_1, x'_k) \\
& = & k+1\,,
\end{eqnarray*}
a contradiction. The other cases similarly lead to contradictions.
\qed

If $|f| = |f'| \le d - 2$, then $\nu({f})=\nu({f'})$ in general no longer implies that $Q_{d}(f)\cong Q_{d}(f')$. For instance, it can be checked that $Q_6(0110) \not\cong Q_6(0100)$ despite the fact that $\nu(0110) = \nu(0100)$. On the other hand we pose the following conjecture.

%For instance, it was checked by computer that $Q_9(0001000)$ is not isomorphic to $Q_9(0001100)$, while $\nu({0001000}) = \nu({0001100}) = 2$. On the other hand we pose:

\begin{conjecture}
\label{conj:nu(f)}
Let $f, f'$ be a non-trivial pair such that $Q_d(f) \cong Q_d(f')$. Then $\nu({f})=\nu({f'})$.
\end{conjecture}

This conjecture has also been verified for all dimensions up to 11 (and all non-trivial pairs $f, f'$); see Appendix~\ref{appendix} for the Sage code. 

\bigskip
\noindent{\bf Acknowledgments:} This research was supported by the bilateral Korean-Slovenian project BI-KR/13-14--005 and the International Research \& Development Program of the National Research Foundation of Korea (NRF) funded by the Ministry of Science, ICT and Future Planning (MSIP) of Korea(Grant number: NRF-2013K1A3A1A15003503). J.A. and S.K. are supported by the Ministry of Science of Slovenia under the grant P1-0297. J.L. and Y.R. are supported by the Basic Science Research Program through the National Research Foundation of Korea funded by the Ministry of Education, Science and Technology grant 2011-00253195.
\ \\\ \\
\noindent J.P.'s research was sponsored by the National Science Foundation under Grant Number DMS-1301692. The United States Government is authorized to reproduce and distribute reprints not-withstanding any copyright notation herein.

\bibliographystyle{acm}

\pagebreak

\appendix
\section{Sage Programs Supporting the Stated Conjectures}
\label{appendix}

In order to test our conjectures we define the function {\em isom\_classes} that returns a dictionary whose entries are nontrivial pairs.
\bigskip\noindent
\begin{python}
def genFib(d,f):

    G = graphs.CubeGraph(d)
    V = [v for v in G.vertices() if v.find(f) != -1]
    G.delete_vertices(V)

    return G

def inv(f):
    return ''.join(['1' if e == '0' else '0' for e in f])

def nu(s):
    cur = s[0]
    ret = 0
    for i in range(1,len(s)):
        if s[i] != cur:
            ret+=1
            cur = s[i]
    return ret

def allStr(k):
    ret = set()
    for f in CartesianProduct(*[['0','1']]*k):
        f = ''.join(f)
        if inv(f) not in ret and f[::-1] not in ret and 
        inv(f[::-1]) not in ret and inv(f)[::-1] not in ret:
            ret.add(f)
    return ret


def isom_classes(d):
    D = {}

    for k in range(3,d):
        for f in allStr(k):
            G = genFib(d,f)
            s = G.canonical_label().graph6_string()
            if s not in D:
                D[s] = [f]
            else:
                D[s]+=[f]
    return D
\end{python}

\pagebreak

\noindent Conjecture~\ref{conj:dim-minus-1} was tested with the following method.

\begin{python}
def testConj1(d):
     D = isom_classes(d)
     for key in D:
         for f1,f2 in Combinations(D[key],2):
             G1 = genFib(d-1,f1)
             G2 = genFib(d-1,f2)
             if not G1.is_isomorphic(G2):
                 return False
     return True

sage: all(testConj1(i) for i in xrange(2,12))
True
\end{python}

\bigskip\noindent
Conjecture~\ref{conj:2/3} was tested with the following method.

\begin{python}
def testConj2(d):
     D  = isom_classes(d)
     for key in D:
         if len(D[key]) > 1 and len(D[key][0]) < (2/3)*(d+1):
             return False
     return True

sage: all([testConj2(i) for i in xrange(3,12)])
True
\end{python}

\bigskip\noindent
Finally, Conjecture~\ref{conj:nu(f)} was tested with the following method.

\begin{python}
def testConj3(d):
    D  = isom_classes(d)
    for key in D:
        if len(set([nu(f) for f in D[key]])) > 1:
            return False
    return True
sage: all([testConj3(i) for i in xrange(3,12)])
True


\end{python}

\end{document}